\documentclass{preprint}

\ifoot{\textbf{Preprint}, Jun 21, 2018}
\titlehead{\textbf{Preprint}, Jun 21, 2018}

\usepackage{float}
\usepackage{amsfonts}
\usepackage{amsmath}
\usepackage{amsthm}
\usepackage{hyperref}

\usepackage{graphicx}
\usepackage{xcolor}  % necessary for *.pdf pics made by fig2eps

\usepackage{booktabs}  % for better style of tables

\usepackage{wrapfig}
\usepackage{amsfonts}
\usepackage{amsmath}

\newtheorem{remark}{Remark}

\newtheorem{form}{Formulation}

\newcommand{\R}{\mathbb R}
\newcommand{\eps}{\varepsilon}
\newcommand{\ep}{\varepsilon}

\newcommand{\bs}[1]{{#1}} % use normal
\newcommand{\bu}{\bs{u}}
\newcommand{\bw}{\bs{w}}

\begin{document}

\title{
  Parallel solution, adaptivity, computational
  convergence, and open-source code of 2d and 3d pressurized
  phase-field fracture problems
}

\author{
Timo Heister
\thanks{Mathematical Sciences, Clemson University,
  Clemson, SC 29634, USA, heister@clemson.edu}
\and
Thomas Wick
\thanks{Institut f\"ur Angewandte Mathematik, Leibniz Universit\"at Hannover,
  Welfengarten 1, 30167 Hannover, Germany,
thomas.wick@ifam.uni-hannover.de}
}

\maketitle

\begin{abstract}
We present a scalable, parallel implementation of a solver for
the solution of a phase-field 
model for quasi-static brittle fracture. The code is available as open source.
Numerical solutions in 2d and 3d with adaptive mesh refinement 
show optimal scaling of the linear solver based on algebraic multigrid,
and convergence of the phase-field model towards exact values of 
functionals of interests such as the crack opening displacement
or the total crack volume.
In contrast to uniform refinement, adaptive mesh refinement allows us to recover optimal convergence rates
for the non-smooth solutions encountered in typical test problems.
We also present numerical studies of the influence of the finite domain size
on functional evaluations used to approximate the infinite domain.

\end{abstract}


%%%%%%%%%%%%%%%%%%%%%%%%%%%%%%%%%%%%%%%%%%%%%%%%%%%%%%%%%%%%%%%%%
\section{Introduction}
In this work, we study the performance of the linear solver and the parallelization of 
our phase-field model for brittle fracture developed in \cite{HeWheWi15} that allows 3d, adaptive solutions
to more than 1024 cores and 100 million degrees of freedom.
The original code was published on \url{https://github.com/tjhei/cracks} 
and is based on the Finite Element library deal.II \cite{dealII85}.
While solver aspects for a similar prolem have been discussed in \cite{FaMau17}, a scalable parallel
implementation is novel.

We consider 2d and 3d benchmark problems in a pressurized fracture setting
(Sneddon's tests \cite{SneddLow69}). 
While seemingly simple, it is a challengening test problem (that is
prototypical for most phase-field fracture configurations) for the following
reasons:
First, the underlying energy functional is non-convex.
Second, the reference problem is defined on an infinite domain, so approximation with a finite domain introduces errors.
Third, due to the smeared interface zone, accurate representations of the total crack volume (TCV) and crack opening displacement (COD)
requires a small phase-field regularization parameter and a very
high resolution mesh around the crack region.

%%%%%%%%%%%%%%%%%%%%%%%%%%%%%%%%%%%%%%%%%%%%%%%%%%%%%%%%%%%%%%%%%
\section{Notation and equations}
Let $B \subset \R^d, d=2,3$ the total domain wherein  
$\mathcal{C}\subset \R^{d-1}$ denotes the fracture
and $\Omega \subset \R^d$ {is} the unbroken domain.
We assume homogeneous Dirichlet conditions on the outer boundary $\partial B$.
Using a phase-field approach, the lower-dimensional fracture $\mathcal{C}$
is approximated by $\Omega_F\subset B$ with the help of an elliptic
(Ambrosio-Tortorelli) functional. 
For fracture formulations posed in a variational setting first proposed in \cite{FraMar98},
fracture regularizations using Ambrosio-Tortorelli functionals were developed in \cite{BourFraMar00}. 
The unknown solution variables 
are vector-valued displacements $\bu:B\to\mathbb{R}^d$ 
and a smoothed scalar-valued indicator phase-field function $\varphi: B\to [0,1]$.
Here $\varphi = 0$ denotes the crack region
and $\varphi = 1$ characterizes the unbroken 
material. The intermediate values constitute a smooth transition zone 
dependent on a regularization parameter $\ep$.
Adding a pressure $p:B\to\mathbb{R}$ 
that acts on 
the fracture boundary was first rigorously modeled and mathematically analyzed in \cite{MiWheWi14_1418}.
Next, the physics require
a crack irreversibility condition (the crack can never heal), which
is an inequality condition in time, i.e, $\varphi \leq \varphi^{old}$,
where $\varphi^{old}$ denotes the previous time step solution.
Consequently, modeling of fracture evolution problems leads to a variational 
inequality system, that is always, due to this constraint,
quasi-stationary or time-dependent.
Let $V:=H^1_0(B),  W:=H^1(B)$ and
$
W_{in}:=\{w\in H^1(B) |\, w\leq \varphi^{old} \leq 1 \text{ a.e. on } B\}
$
be the function spaces to state the variational formulation.
In the following, we denote the $L^2$ scalar product with $(\cdot, \cdot)$.
The Euler-Lagrange system for pressurized phase-field fracture reads \cite{HeWheWi15,MiWheWi14_1418}:
\begin{form}
\label{form_1}
Let $p\in L^{\infty}(B)$ be given. For the loading steps $n=1,2,3,\ldots$: Find vector-valued displacements and a
scalar-valued phase-field variable $\{\bu,\varphi\} :=\{\bu^n,\varphi^n\} \in \{u_D + V\} \times W$ such that
\begin{align}
  &\Bigl(\big( (1-\kappa) \tilde{\varphi}^2  +\kappa \big)\;\sigma(\bu), e( {\bw}
  )\Bigr)  
+({\varphi}^{2} p, \mbox{div }  {\bw}) 
=0, \label{E22ATW}\\ 
&(1-\kappa) ({\varphi} \;\sigma(\bu):e( \bu), \psi {-\varphi}) 
+  2 ({\varphi}\;  p\; \mbox{div }  \bu,\psi{-\varphi})
+  G_c  \Bigl( -\frac{1}{\ep} (1-\varphi,\psi{-\varphi}) + \ep (\nabla
\varphi, \nabla (\psi - {\varphi}))   \Bigr)  \geq  0,\label{E33ATW}
\end{align}
for all $\bw\in V$ and $\psi \in W_{in}\cap L^{\infty}(B)$.
Here, $G_c$ is the critical energy release rate, and
we use the well-known law for the linear stress-strain relationship
$
\sigma := \sigma ( {u}) = 2\mu_s e( {u}) + \lambda_s \text{tr}e( {u})I,
$
where $\mu_s>0$ and $\lambda_s>0$ denote the Lam\'e coefficients, 
$e( {u}) = \frac{1}{2}(\nabla {u} + \nabla {u}^T)$ 
is the linearized strain tensor 
and $I$ is the identity matrix, $\kappa > 0$. Finally, $\tilde\varphi$ is 
a linear extrapolation in time developed in \cite{HeWheWi15} in order 
to convexify the above problem. 
\end{form}

%%%%%%%%%%%%%%%%%%%%%%%%%%%%%%%%%%%%%%%%%%%%%%%%%%%%%%%%%%%%%%%%%
\section{Solution algorithms and parallel framework}
To solve Formulation \ref{form_1}, we employ a semi-smooth Newton method that was developed for 
phase-field fracture in \cite{HeWheWi15} and combines two 
Newton methods: solving the nonlinear problem and treating 
the irreversibility constraint. 
The spatial discretization is based 
on a Galerkin finite element 
scheme, introducing $H^1$ conforming discrete spaces $V_h\subset V$ and 
$W_h \subset W$ consisting of bilinear/trilinear functions $Q_1^c$ on
quadrilaterals / hexahedra, respectively.
The discretization parameter is denoted by $h$.
The code is fully parallelized using MPI by building on the deal.II finite 
element library \cite{dealII85}. The adaptive meshes are handled by p4est \cite{p4est}
and the linear algebra is built on Trilinos \cite{Trilinos-Overview}. 
This {software framework} is discussed in \cite{BangerthBursteddeHeisterKronbichler11}.

From our original work \cite{HeWheWi15}, we extended the active set strategy with a method to detect and constrain
alternating active set indices to avoid cycles of the method similar to \cite{CurtisActiveSet}.
Additionally, we extended our refinement strategy from not only refining in
the crack region to enforce $\eps = 2h$
to also include a gradient jump estimator for the displacements. This is
required to achieve rigorous computational convergence of
the benchmark problem discussed here.

%\subsection{A combined Newton method: solving the nonliear system and the
%  crack irreversibility constraint}
%The algorithm reads:
%\begin{Algorithm}
%\label{final_algo}
%Repeat for $k=0, \dots$ until the active set $\mathcal{A}_k$ does not 
%change and $ \widetilde{R}(U_h^k) < \operatorname{TOL}$:
%\begin{enumerate}
% \item Assemble residual $R(U_h^k)$
% \item Compute active set $\mathcal{A}_k = \{i \mid  (B^{-1})_{ii} (R_k)_i + c (\delta U_h^k)_i > \eps \}$
% \item Assemble matrix $G = \nabla^2 E_\ep(U_h^k)$ and right-hand side $F = -\nabla E_\ep(U_h^k)$
% \item Eliminate rows and columns in $\mathcal{A}_k$ from $G$ and $F$ to obtain $\widetilde{G}$ and $\widetilde{F}$
% \item Solve linear system with a block-preconditioned GMRES scheme: $\widetilde{G} \delta U_k = \widetilde{F}$, i.e, find $\delta U_h^k \in V_h \times W_h$ with
%\begin{equation}
%\label{linear_system_in_Newton}
%  \nabla^2 E_\ep( {U_h^k})(\delta U_h^k, {\Psi})  = - \nabla E_\ep( {U_h^k})( {\Psi}) \quad\forall\Psi\in V_h\times W_h,
%\end{equation}
%where $\nabla^2 E_\ep$ and $\nabla E_\ep$ 
%are the Hessian matrix and gradient respectively.
%
% \item Find a step size $0 < \omega \leq 1$ using line search to get
% \[ 
%  U_h^{k+1} = U_h^k + \omega \delta U_h^k,
% \]
% with $ \widetilde{R}(U_h^{k+1})<\widetilde{R}(U_h^k)$.
%\end{enumerate}
%\end{Algorithm}

\subsection{Linear iterative solution}
The linear systems arising at each Newton step 
are solved iteratively using a GMRES scheme
with a block diagonal preconditioner $P^{-1}$:
\begin{equation*}
 \begin{pmatrix}
  M_{uu} & 0 \\ M_{\varphi u} & M_{\varphi\varphi}
 \end{pmatrix}
 \begin{pmatrix}
 \delta \bu \\ \delta \varphi
 \end{pmatrix}
 =
 \begin{pmatrix}
  F_u \\ F_\varphi
 \end{pmatrix}
\quad\text{with} \quad
 P^{-1}=\begin{pmatrix}
  \tilde{M}_{uu}^{-1} & 0 \\ 
  0 & \tilde{M}_{\varphi\varphi}^{-1}
 \end{pmatrix}
.
\end{equation*}
Using the basis $\{\psi_i| i=1,\ldots,N\}$ in $V_h\times W_h$ with 
$\psi_i = (\chi_i^u,0)^T,i=1,\ldots,N_u$ 
and $\psi_{N_u+i} = (0, \chi_i^{\varphi})^T,i=1,\ldots,N_{\varphi}$ and 
$N = N_u + N_{\varphi}$, we have
specifically the entries: 
\begin{align*}
 (M_{uu})_{i,j}
&=
\Bigl(\big( (1-\kappa) {\tilde\varphi}^2  +\kappa \big)\;\sigma({\mathbf{\chi}}_j^u), e( {\mathbf{\chi}}_i^u )\Bigr)
+ (\sigma({\mathbf{\chi}}_j^u), e( {\mathbf{\chi}}_i^u )), \\
%%%
(M_{\varphi u})_{i,j}
&= 
2(1-\kappa) ({\varphi} \; \sigma({\mathbf{\chi}}_j^u):e( u),\chi_i^{\varphi})
-  2(\alpha -1) p ({\varphi}\; \mbox{div} ( {\mathbf{\chi}}_j^u),\chi_i^{\varphi}),  \\
%%%
%M^{u\varphi}_{i,j}
%&=0, \\
%%%
(M_{\varphi\varphi})_{i,j}
&= (1-\kappa) ( \sigma(u):e(u) \chi_j^{\varphi},\chi_i^{\varphi})
 -  2(\alpha -1) p (\mbox{div}( \bu) \chi_j^{\varphi},\chi_i^{\varphi})
+  G_c  \Bigl( \frac{1}{\varepsilon} (\chi_j^{\varphi},\chi_i^{\varphi}) + \varepsilon (\nabla
\chi_j^{\varphi}, \nabla \chi_i^{\varphi})   \Bigr). 
%%%%
\end{align*}
The block $M_{u\varphi}$ is zero due 
$\tilde\varphi$ in Equation \eqref{E22ATW}.

We assume the existence of spectrally equivalent approximations $\tilde{M}_{uu}^{-1}$ and $\tilde{M}_{\varphi\varphi}^{-1}$,
which correspond to linear elasticity and a mixture of a Laplacian and mass matrix, respectively. The eigenvalues are given by
the generalized eigenvalues of the systems $\tilde{M}_{kk} = \lambda {M}_{kk}$ ($k=u, \varphi$). We are approximating
the blocks using a single V-cycle of algebraic multigrid (by Trilinos ML, 
\cite{Trilinos-Overview}).
Assuming the multigrid gives approximations independent of the mesh size $h$,
the whole solver scheme is nearly optimal and independent
of the number of processors and mesh size, as can be seen in Table~\ref{fig:solver}.

\renewcommand{\arraystretch}{1.3}

\begin{table}[h!]
\begin{center}
{%\tt   % see code/solver_it/  
\begin{tabular}{|rr|rrrrrrrr|} \hline
  & & &  &  &  & NP &  & & \\
  ref & Dofs & 16 & 32 & 64 & 128 & 256 & 512 & 1024 & 2048 \\ \hline
%   1 & 867 & 2 & 2 & 2 & 2 & 2 & 2 & 2 & 2 \\
%   2 & 3267 & 10 & 10 & 9 & 9 & 8 & 8 & 8 & 8 \\
%   3 & 12675 & 13 & 13 & 13 & 13 & 12 & 12 & 19 & 19 \\
%   4 & 49923 & 17 & 17 & 16 & 17 & 14 & 15 & 14 & 14 \\
  5 & 198'147 & 19 & 19 & 19 & 18 & 18 & 18 & 17 & 19 \\
  6 & 789'507 & 24 & 23 & 25 & 24 & 25 & 23 & 24 & 23 \\
  7 & 3'151'875 & 33 & 27 & 29 & 28 & 27 & 27 & 25 & 37 \\
  8 & 12'595'203 & - & - & 31 & 31 & 33 & 32 & 30 & 32 \\
  9 & 50'356'227 & - & - & 43 & 43 & 44 & 48 & 40 & 52 \\ \hline
  \end{tabular}
} 
\end{center}
\caption{Number of GMRES iterations of a single Newton step for the Sneddon 2d test with global refinement.
Iterations are nearly independent of problem size ($h$) and number of
processors $NP$. The relative residual is 1e-8.}
\label{fig:solver}
\end{table}

\renewcommand{\arraystretch}{1.0}

\begin{remark}
The parallelized linear solution of phase-field fracture was already
implemented in \cite{HeWheWi15} for 2d problems, 
but a computational analysis of the
performance and the extension to 3d settings was missing up to now in the existing literature.
\end{remark}

\newpage
%%%%%%%%%%%%%%%%%%%%%%%%%%%%%%%%%%%%%%%%%%%%%%%%%%%%%%%%%%%%%%%%%
\section{Numerical tests: Sneddon 2d and 3d}
We perform numerical tests in 2d and 3d.
The benchmarks are based on 
the theoretical calculations of \cite{SneddLow69}[Section 2.4 and Section 3.3].
A (constant) pressure $p = 10^{-3}Pa$ 
causes 
the fracture to change its width but not the length to form a penny-shaped crack.
The initial crack of length $l_0=1.0$ is described with the help of the 
phase-field function $\varphi$. The domains are varied as $(-K,K)^d, d=2,3, K=5,10,20,40$ to
approximate the infinite domain of the benchmark. We choose $G_c = 1$, Youngs' modulus as
$E=1$ and Poisson's ratio is $\nu = 0.2$. The regularization parameters are
chosen as $\eps = 2h$ and $\kappa = 10^{-12}h$. 
For more details of our computational configuration we refer 
the reader to \cite{HeWheWi15}[Section 5.3] (2d) and 
\cite{WheWiWo14}[Section 5.4] (3d).

We study the total crack volume (TCV) with analytical solutions from \cite{SneddLow69}, respectively:
{\small 
\begin{equation}
TCV_{h,\epsilon} = \int_\Omega u \cdot \nabla \phi dx, \;
 TCV_{2d} = \int_x 2u_y(x) dx = \frac{2\pi p l_0^2(1-\nu^2)}{E}, \;
 TCV_{3d} = \int_x\int_y 2u_z(x,y) dxdy = \frac{16 p l_0^3(1-\nu^2)}{3E}
\end{equation}
}
where $l_0=1.0$, $E=1.0$, $p=1e-3$, $\nu=0.2$.
We note that the crack opening displacement (COD) is computed by
\[
 u_n(v) = \frac{cp l_0 (1-\nu^2)}{E}\sqrt{1-\left(\frac{\rho}{l_0}\right)^2},
 \quad\text{with the radius } \; \rho=\|x\|_{l^2}, x\in R^d 
\]
with $n=y, v=x, c = 2$ in $2d$ and $n=z, v=(x,y), c = 4/\pi$ in $3d$.

\subsection{Convergence of total crack volume (TCV) in 2d}
\label{sec_41}

\begin{figure}[h!]
\centering
\includegraphics[width=8cm]{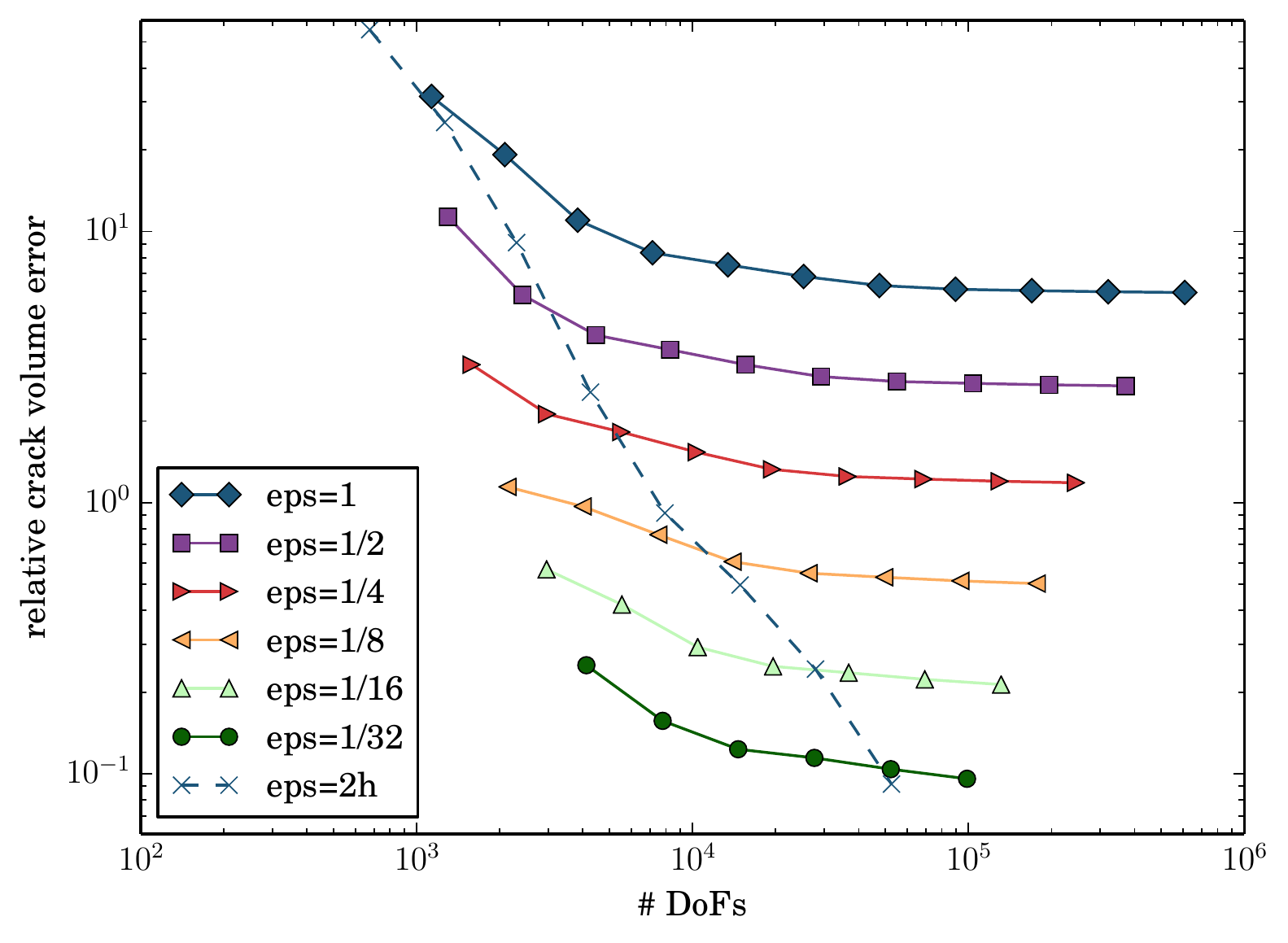}
\includegraphics[width=8cm]{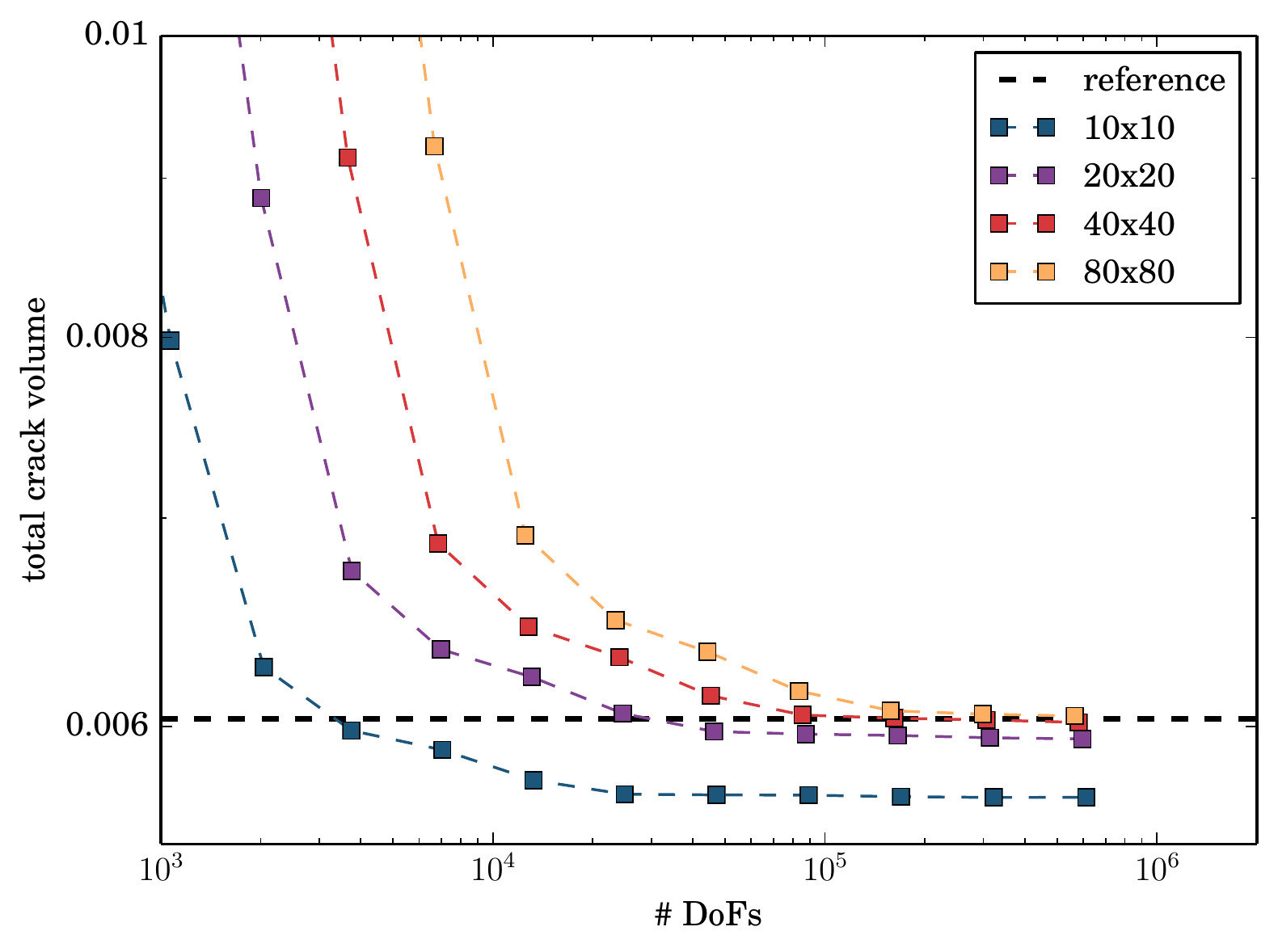}
\caption{Left (Section \ref{sec_41}): The convergence of $\ep$ on the 40x40 domain is computationally
analyzed.
Here we observe linear convergence in epsilon.
Right (Section \ref{sec_42}): dependence of functional value evaluations on the domain
  size.}
\label{fig:sneddon-2d-tcv}
\end{figure}

The first set of computations computes the error of the TCV in 2d for fixed values of $\eps$ while
adaptively refining the solutions. Convergence for $h\rightarrow 0$ and $\eps \rightarrow 0$ are
clearly visible; see Fig.~\ref{fig:sneddon-2d-tcv}, left.
This computation demonstrates: First, the phase-field regularization is a valid and convergent approximation of the fracture problem.
Second, our approach for adaptive refinement is effective and the choice of $\eps=2h$ combined with adaptivity is
a valid implementation strategy, that gives accurate solutions with a few number of unknowns.

\newpage
\subsection{Influence of finite domain size in 2d}
\label{sec_42}
The benchmark problem is stated for an infinite domain, so convergence is only obtained
until this part of the error becomes dominant as can be seen in Fig.~\ref{fig:sneddon-2d-tcv} (right).
A larger size of the domain would raise the cost of the discretization immensly without adaptive 
mesh refinement, while we obtain more accurate solutions on a larger domain with about $10^5$ unknowns
and our refinement strategy.
The extrapolated value of the TCV compared to the exact value of 6.0319E-03 
has an error of 5.6\%, 1.5\%, 0.5\%, and 0.1\% for the domain size of 10x10, 20x20, 40x40, and 80x80, respectively.

% \begin{wrapfigure}{R}{0.4\textwidth}
% {\tt
% \footnotesize
% \begin{tabular}{r|r|r} % data taken from paper_2/sneddon_2d_domain/sneddon_cod.ods
%   Domain & TCV & rel. Error \\ \hline
%   10x10 & 5.6942E-03 & 5.60\% \\
%   20x20 & 5.9393E-03 & 1.53\% \\
%   40x40 & 6.0014E-03 & 0.51\% \\
%   80x80 & 6.0383E-03 & 0.11\% \\
%   reference & 6.0319E-03 &  \\
%   \end{tabular}
% } 
% \caption{Sneddon 2d, influence of domain size}
% \label{fig:domain-size}
% \vspace{-2em}
% \end{wrapfigure}

%Same data, but the value given in the following table is extrapolated from the three finest computations in Fig.~\ref{fig:domain-size}

\subsection{Convergence of the 3d benchmark}
\label{sec_43}
\begin{figure}[h!]
\begin{minipage}[t]{0.38\textwidth}
\vspace{1em}
%\fbox
{
 \includegraphics[width=\textwidth]{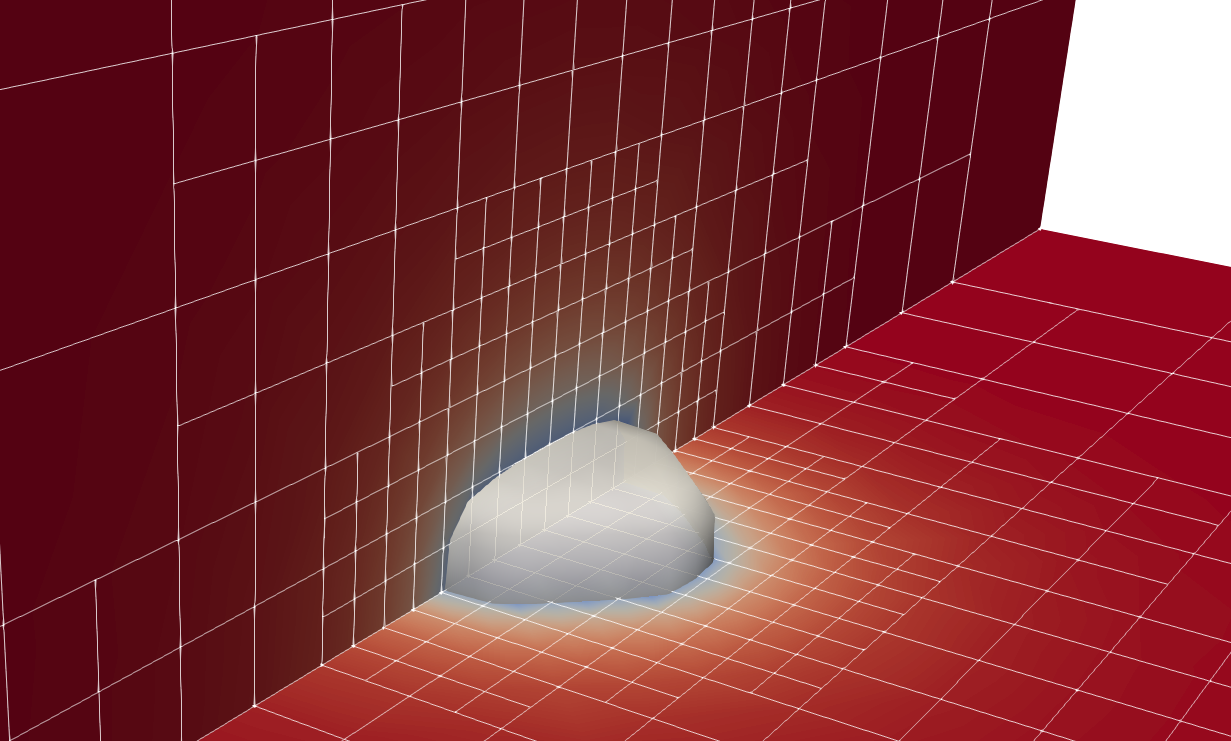}}
\end{minipage} 
\hspace{1em}
\begin{minipage}[t]{0.6\textwidth}
\vspace{1em}
%\fbox
{%\tt % see paper_2/sneddon_3d/3d.prm and output cracks.o9389527 from palmetto
\footnotesize
\begin{tabular}{|rrrrrr|} \hline
  ref & $10l_0/h$ & \# DoFs & eps & TCV & Error \\ \hline
%  1 & 2 & 500 & 1.7321E+01 & 7.4519E-02 & 1355.44\% \\
%  2 & 4 & 1,220 & 8.6603E+00 & 5.3907E-02 & 952.86\% \\
%  3 & 8 & 4,340 & 4.3301E+00 & 1.4066E-02 & 174.72\% \\
  4 & 16 & 12,436 & 2.1651E+00 & 3.5158E-02 & 586.67\% \\
  5 & 32 & 35,940 & 1.0825E+00 & 2.0078E-02 & 292.14\% \\
  6 & 64 & 105,740 & 5.4127E-01 & 1.0149E-02 & 98.22\% \\
  7 & 128 & 313,668 & 2.7063E-01 & 6.9110E-03 & 34.98\% \\
  8 & 256 & 946,852 & 1.3532E-01 & 5.8082E-03 & 13.44\% \\
  9 & 512 & 2,938,012 & 6.7658E-02 & 5.3764E-03 & 5.01\% \\
  10 & 1024 & 9,318,916 & 3.3829E-02 & 5.2131E-03 & 1.82\% \\
  11 & 2048 & 30,330,756 & 1.6915E-02 & 5.1567E-03 & 0.72\% \\
  12 & 4096 & 100,459,828 & 8.4573E-03 & 5.1352E-03 & 0.30\% \\ \hline
  & \multicolumn{3}{l}{reference (Sneddon)}    & 5.1200E-03 &  \\ \hline
  \end{tabular}
}
\end{minipage} 

\caption{Section \ref{sec_43}: Sneddon 3d adaptive convergence for a 10x10x10 domain, $l_0=1$, using $\ep=2h$. Left: solution with mesh and
isosurface for $\varphi=0.3$. Right: convergence table.}
\label{fig:sneddon3d-results}
\end{figure}

We now concentrate on the convergence of the TCV error for the 3d version of the same benchmark problem.
Even with adaptive refinement, an accurate solution requires in the order of $100$ million DoFs 
and substantial computing power (here run on $1\, 024$ cores); see Figure~\ref{fig:sneddon3d-results}.
On the one hand, this is a well-known challenge for smeared interface
approaches. On the other hand, 3d computations for 
moving interface problems are still a challenge for `exact' interface 
representations such as extended/generalized finite elements, cut cell methods, etc..
Therefore, these results show that phase-field approaches can give comparable results to other
numerical methods.

\renewcommand{\arraystretch}{1.3}

\renewcommand{\arraystretch}{1.0}

%\newpage
\subsection{Adaptive convergence of crack opening displacement (COD)}
\label{sec_44}

\begin{figure}[h!]
\begin{center}
\includegraphics[width=0.65\textwidth]{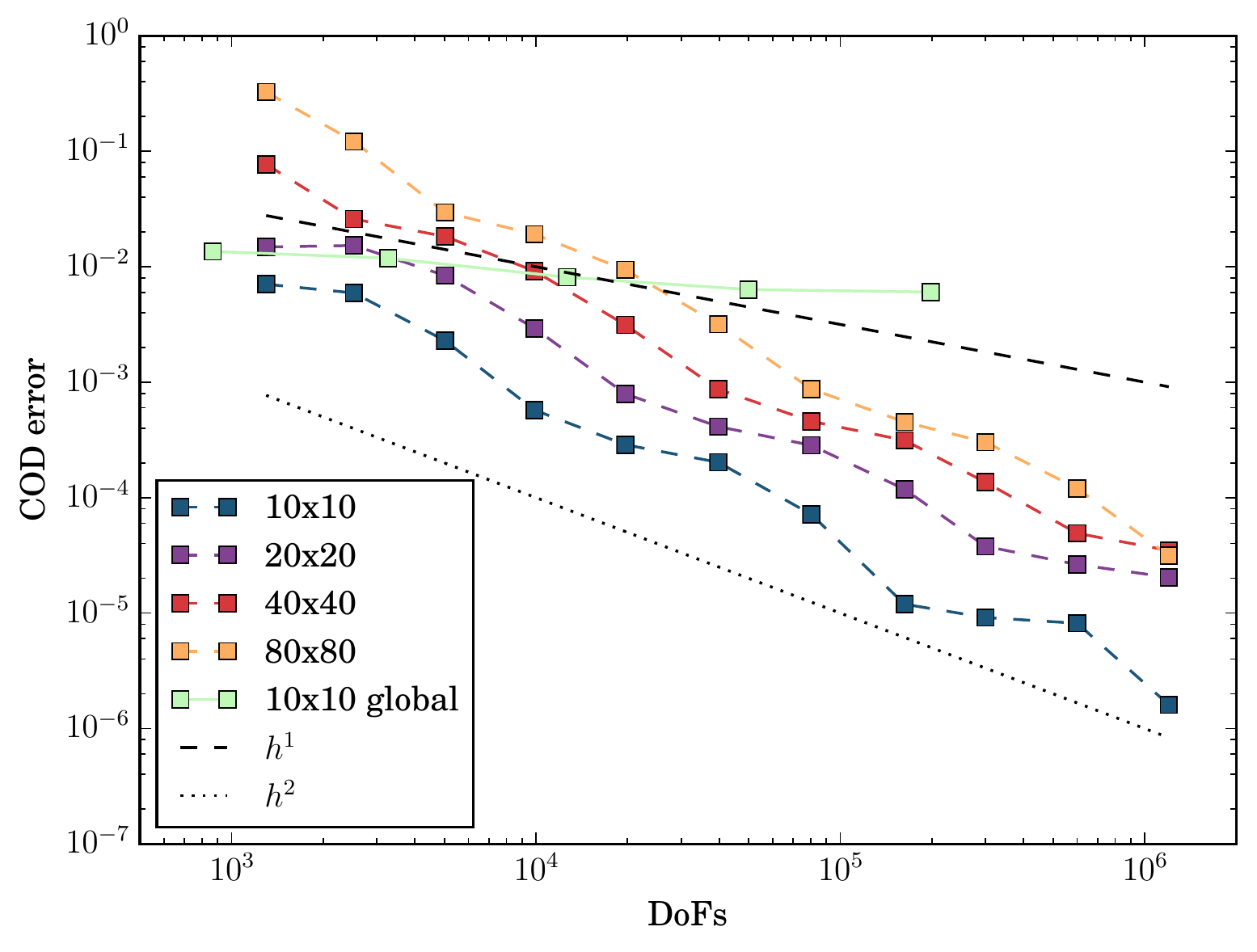} 
\end{center}
\caption{Section \ref{sec_44}: Sneddon 2d error convergence to extrapolated COD values for a fixed domain size. We recover quadratic convergence in $h$ using our adaptive refinement strategy.
Errors are increasing for a fixed number of unknowns for a larger domain size, as the crack is resolved with fewer cells.
}
\label{fig:sneddon-2d-error}
%\end{wrapfigure}
\end{figure}

Finally, we study the crack opening displacement (COD) in 2d to compare adaptive vs. global refinement. 
For each domain size, the reference value is obtained using Richardson extrapolation of the numerical solution.
As the solution
is discontinuous, we see suboptimal convergence of the solution using global refinement, while we recover optimal
rates using our adaptive scheme; see Figure~\ref{fig:sneddon-2d-error}.

Our results suggest
a convergent method and lead us to expect an error estimate of the kind (similar to \cite{Wi16_dwr_pff})
\begin{align*}
 \| u_{h,\ep} - u_{ref} \| &\leq \| u_{h,\ep} - u_{ref,\ep} \| + \| u_{ref,\ep} - u_{ref} \| \\
 &\leq C_1 \inf_v \| v - u_{ref,\ep} \| + C_2\ep, \quad C_1,C_2 > 0,
\end{align*}
where the error of the discrete solution $u_{h,\ep}$ to the exact solution $u_{ref}$ is given by the sum of best-approximation error
of the finite element method and an error term (i.e., a model error) introduced 
by the $\ep$ regularization converging linearly in $\ep$.

%\begin{wrapfigure}{R}{0.4\textwidth}

%%%%%%%%%%%%%%%%%%%%%%%%%%%%%%%%%%%%%%%%%%%%%%%%%%%%%%%%%%%%%%%%%
\section*{Acknowledgements}
  This work is supported by the 
German Priority Programme 1748 (DFG SPP 1748) 
Reliable Simulation Techniques in Solid Mechanics. Development of Non-standard Discretization Methods, Mechanical and Mathematical Analysis.
Timo Heister was partially supported by the Computational Infrastructure in
Geodynamics initiative (CIG), through the NSF under Award EAR-0949446 and The
University of California – Davis, by the NSF Award DMS-1522191, and
by Technical Data Analysis, Inc through US Navy SBIR N16A-T003.
Clemson University is acknowledged for generous allotment of compute time on Palmetto cluster.

\newpage
%%%%%%%%%%%%%%%%%%%%%%%%%%%%%%%%%%%%%%%%%%%%%%%%%%%%%%%%%%%%%%%%%

\end{document}